\documentclass[12pt]{article}
\title{Tracing Pareto-optimal points for multi-objective shape optimization applied to electric machines}

\usepackage{a4,epsfig}
\usepackage{amsfonts,amsmath}
\usepackage{subcaption}
\usepackage{algorithm}
\usepackage{algpseudocode}
\usepackage{pgfplotstable}

\newtheorem{remark}{Remark}[section]

\newcommand{\R}{\mathbb{R}}

\newcommand{\black}{\color{black}}

\numberwithin{equation}{section} 

\begin{document}

\setcounter{page}{1}

\title{Tracing Pareto-optimal points for multi-objective shape optimization applied to electric machines}

\author{Alessio Cesarano$^1$, Peter Gangl$^1$}
\date{$^1$Johann Radon Institute for Computational and
  Applied Mathematics, \\
Altenberger Stra{\ss}e 69, 4040 Linz, Austria }

\maketitle

\begin{abstract}
In the context of the optimization of rotating electric machines, many different objective functions are of interest and considering this during the optimization is of crucial importance. While evolutionary algorithms can provide a Pareto front straightforwardly and are widely used in this context, derivative-based optimization algorithms can be computationally more efficient. In this case, a Pareto front can be obtained by performing several optimization runs with different weights. In this work, we focus on a free-form shape optimization approach allowing for arbitrary motor geometries. In particular, we propose a way to efficiently obtain Pareto-optimal points by moving along to the Pareto front exploiting a homotopy method based on second order shape derivatives.
\end{abstract}

\section{Introduction}\label{sec:intro}

Free-form shape optimization \cite{sozo} allows for arbitrarily shaped geometries that have the same topology as the initial design and, thus, yields a larger design space compared to geometric parameter optimization. In many engineering applications, one wants to optimize multiple, often conflicting, cost functions and thus is interested in the set of (locally) Pareto-optimal designs. Within the class of evolutionary algorithm \cite{deb}, it is straightforward to extract Pareto fronts by selecting those designs which are not dominated by any of the other evaluated designs. While derivative-based algorithms can yield optimized designs for single-objective optimization problems much faster, obtaining Pareto fronts in the multi-objective setting requires dedicated strategies such as the weighted sum approach \cite{das} or approaches based \black{on} descent directions with respect to multiple objective functions \cite{gangl}. These approaches, however, still require the solution of a full single-objective optimization for each point on the Pareto front. In order to overcome this efficiency issue, methods for tracing along the Pareto front have been introduced in \cite{schütze, hillermeier, schmidt, doganay}, i.e., to parametrize the Pareto front and sequentially \black{compute} Pareto-optimal points without running a new optimization. In contrast to the widely used weighted sum approach, these methods also allow for finding points on non-convex regions of the Pareto front and to have better control over the spacing of the points on the Pareto front.

In this work, we extend the Pareto tracing methods for efficiently obtaining many points on the Pareto front to the case of free-form shape optimization under PDE constraints. In particular, we apply the methodology to the shape optimization of a synchronous reluctance machine.

The rest of this paper is organized as follows: In Sect.~\ref{cesarano:sec_probDesc} we describe the physical model we consider for our numerical experiments. We shortly introduce the mathematical tools used for this work, i.e. the concept of homotopy methods and shape Newton methods, in Sects.~\ref{cesarano:sec_homotopy} and \ref{cesarano:sec_shapeNewton}, respectively, before presenting their application to the physical problem and discussing our numerical results in Sect.~\ref{cesarano:sec_paretomotor}.

\section{Problem description}\label{cesarano:sec_probDesc}

The machine under study is used in X-ray generation machines and it has a large air gap between rotor and stator. In reality the rotor of the machine is in vacuum and is separated from the stator by glass, which has magnetic properties of air. We choose to restrict the shape deformation only to the rotor, changing its distribution of ferromagnetic material and air. We intend to, at the same time, minimize both the negative torque (i.e. maximize the torque) and the weight of the rotor, i.e., we consider the bi-objective optimization problem

\begin{equation}\label{cesarano:eq_arkkio}
            \underset{\Omega}{\mbox{min }} \left(\mathcal J_1(\Omega), \mathcal  J_2(\Omega) \right)^\top = \left( - \nu_0 \: L \frac{1}{r_2-r_1}\int_{\Omega_g} (\nabla u_\Omega )^\top Q \nabla u_\Omega \; \mbox dx , \text{Vol}(\Omega) \right)^\top\\
 \end{equation}
 \begin{equation}
\label{cesarano:eq_magnetostatics}
	\mbox{where } u_\Omega \in H^1_0(D): \int_{D} \nu_{\Omega}(x,|\nabla u_\Omega|)\nabla u_\Omega \cdot \nabla v \, \mbox dx = \int_D J_z v \, \mbox dx \; \forall v \in H^1_0(D).
\end{equation}
Here, \eqref{cesarano:eq_magnetostatics} represents the equation of 2D nonlinear magnetostatics.
In \eqref{cesarano:eq_arkkio} we use the Arkkio formula in order to calculate the torque. We denote as $u$ the $z$-component of the magnetic vector potential, being
\begin{equation}\label{cesarano:eq_acurl}
B = \mbox{curl}  A = \mbox{curl} (0, 0, u(x_1, x_2))^\top.
\end{equation}

In Arkkio's formula, the dot product between the tangential and normal component of $B$ is needed. The $2 \times 2$ matrix $Q$
\begin{equation}\label{cesarano:eq_qu}
   Q(x_1,x_2) = \frac{1}{\sqrt{x_1^2+x_2^2}} \begin{pmatrix}x_1 x_2 & \frac{x_2^2-x_1^2}{2} \\
        \frac{x_2^2-x_1^2}{2} & - x_1 x_2 \end{pmatrix}
\end{equation}
encodes this dot product, taking into account the transformation from the scalar $u$ to the $2$-dimensional $B$ field. Moreover, $L$ is the length of the motor in $z$ dimension, $\Omega_g$ is a ring in the air gap with inner and outer radii $r_1<r_2$, respectively.

The magnetic reluctivity $\nu$ is defined piecewise. It is constant in air while it depends on the $|B|=|\nabla u|$ in the ferromagnetic material,
\begin{equation}
	\nu_{\Omega}(x,|\nabla u|) := 
	\begin{cases}
		\hat \nu(|\nabla u|), & x \in \Omega, \\
        \nu_0, & x \in  D \setminus  \overline \Omega .
    \end{cases}
\end{equation}
Here we denote with $D$ the whole domain and with $\Omega$ the ferromagnetic subdomain. The current density $J_z$ is piecewise costant and non-zero only in the coils. The machine has 1 pole pair and the angle between the fields deriving from the currents and the $y$-axis is 45 degrees. The maximum amplitude is 10 Ampere.

\section{Homotopy methods} \label{cesarano:sec_homotopy}

Homotopy methods have been introduced for solving (systems) of nonlinear equations by connecting the problem of interest smoothly to a problem whose solution is known or readily computed and following the solution path \cite{allgower}. More precisely, let $F : \R^n \rightarrow \R^n$ and $x^{(0)} \in \R^n$ be given. We look for a root satisfying
the equation
\begin{equation} \label{cesarano:eq_Fxzero}
    F(x) = 0.
\end{equation}
Often, a solution to problem \eqref{cesarano:eq_Fxzero} cannot be found from the given initial guess $x^{(0)}$. A homotopy map $H: \R^n \times \R \rightarrow \R^n$ defines a family of equations

\begin{equation}\label{cesarano:eq_Hxzero}
    H(x,t) = 0
\end{equation}
parametrized by the homotopy parameter $t \in [0,1]$. The homotopy map should be chosen in such a way that the initial problem $H(x,0) = 0$ has a known solution $x^{(0)}$ or its solution can be easily computed from the initial guess $x^{(0)}$, and the equation $H(x,1) = 0$ is equivalent to the original equation \eqref{cesarano:eq_Fxzero}.
A popular example is the convex homotopy 
\begin{equation} \label{cesarano:eq_convHomo}
	H(x,t) = tF(x) + (1-t) G(x)
\end{equation}
where $G:\R^n \rightarrow \R^n$ and $x^{(0)}$ with $G(x^{(0)}) = 0$ are given. The general idea of homotopy methods then is, starting out from $t=0$, to follow the solution path, most often using a predictor-corrector method. {\black The concept of homotopy methods was extended to optimization problems in \cite{dunlavy}.}

We want to use this approach in the context of shape optimization, meaning that the functions $G$ and $F$ will be replaced by optimality conditions for a simple and desired optimization problem, respectively.
For the predictor step, first or higher order shape derivatives may be used. For the corrector step, it makes sense to use the contraction properties of Newton-like second order methods in order to be sure to be and remain close enough to the homotopy curve. In our framework, this lead us to the use of a shape Newton method.

\section{Shape Newton} \label{cesarano:sec_shapeNewton}

Given a shape function $\mathcal J=\mathcal J(\Omega)$, smooth vector fields $V, W : \R^d \rightarrow \R^d$, parameters $s, r \geq 0$ and corresponding perturbed domains $\Omega_{s}:= (\text{id}+s V)(\Omega)$ and $\Omega_{s,r}:= (\text{id}+s V + r W)(\Omega)$, we will make use of the first and second order shape derivatives defined as
\begin{equation} \label{cesarano:eq_def_d2J}
d \mathcal J(\Omega)(V) := \left.\left(\frac{d}{ds} \mathcal J(\Omega_{s}) \right) \right\rvert_{s=0}, \;
    d^2 \mathcal J(\Omega)(V, W) := \left.\left(\frac{d^2}{dsdr} \mathcal J(\Omega_{s, r}) \right) \right\rvert_{s=0, r=0}.
\end{equation}
The main idea of most first order shape optimization algorithms consists in repeatedly finding a descent direction $V$ such that $d\mathcal J(\Omega)(V)<0$ and updating the current domain by the action of this vector field to get $\tilde{\Omega} = (\text{id} + s V)(\Omega)$ as the new iterate. A descent direction can be obtained by solving an auxiliary boundary value problem to find $V \in H$ satisfying
\begin{equation} \label{cesarano:eq_desc_direction}
    b(V, W) = - d \mathcal J(\Omega)(W)
\end{equation}
for all $W \in H$ for a given vector-valued Hilbert space $H$ and a positive definite bilinear form $b(\cdot, \cdot) : H \times H \rightarrow \R$. Replacing the bilinear form $b(\cdot, \cdot)$ in \eqref{cesarano:eq_desc_direction} by the (symmetric) second order shape derivative defined in \eqref{cesarano:eq_def_d2J} for the problem at hand and the current domain $\Omega$ yields the shape-Newton system to find $V \in H$ satisfying
\begin{equation}\label{cesarano:eq_shapeNewton}
    d^2\mathcal J(\Omega)(V, W) = - d \mathcal J(\Omega)(W)
\end{equation}
for all $W \in H$. The operator on the left hand side of \eqref{cesarano:eq_shapeNewton}, however, is not invertible. In fact, from Hadamard's structure theorem \cite{sozo} we know that only normal deformations at the boundary lead to changes in the cost function. Any other deformations, such as interior or tangential deformations, are in the kernel of the operator. 

After finite element discretization, we remove interior deformations by restricting \eqref{cesarano:eq_shapeNewton} to the boundary $\partial \Omega$ and we filter out tangential movements by imposing the pointwise constraints $(V \cdot \tau)(x_k) = 0$ for all mesh vertices $x_k$ on $\partial \Omega$ where $\tau$ denotes the tangential vector.

\section{Application to model problem} \label{cesarano:sec_paretomotor}
Here, we combine the methods introduced in Sects.~\ref{cesarano:sec_homotopy} and \ref{cesarano:sec_shapeNewton} and apply them to the bi-objective PDE-constrained shape optimization problem \eqref{cesarano:eq_arkkio}--\eqref{cesarano:eq_magnetostatics}.

\subsection{Pareto tracing for shape optimization} \label{cesarano:sec_paretoshape}
As outlined in Sect.~\ref{cesarano:sec_homotopy}, we define the convex homotopy between the optimality conditions for cost functions $\mathcal J_1(\Omega)$ and $\mathcal J_2(\Omega)$ defined in \eqref{cesarano:eq_arkkio},
\begin{equation} \label{cesarano:eq_H_motor}
    \mathcal H(\Omega, t)(\cdot) := (1-t) d \mathcal J_1(\Omega)(\cdot) + t d \mathcal J_2(\Omega)(\cdot),
\end{equation}
where $(\cdot)$ represents some vector-valued test functions. For the shape derivatives of the torque and the volume, see e.g. \cite{gangl}. Given a solution of $\mathcal H(\Omega, 0)(V)=0$ for all test functions $V$, i.e., a locally optimal design with respect to the torque, we adaptively choose a sequence of intermediate homotopy values $0<t_k \leq 1$ and obtain the next Pareto optimal design $\Omega^{(k+1)}$ as the solution of $\mathcal H(\Omega, t_{k+1})(\cdot)=0$ by a shape Newton method starting out from $\Omega = \Omega^{(k)}$, i.e., by repeatedly solving systems of the form \eqref{cesarano:eq_shapeNewton} as described in Sect.~\ref{cesarano:sec_shapeNewton}. Note that the procedure can also be extended to more than two cost functions, see \cite{hillermeier}.

\subsection{Numerical results} \label{cesarano:sec_numerics}
\begin{figure}
\begin{tabular}{cc}
    \includegraphics[width=0.35\textwidth]{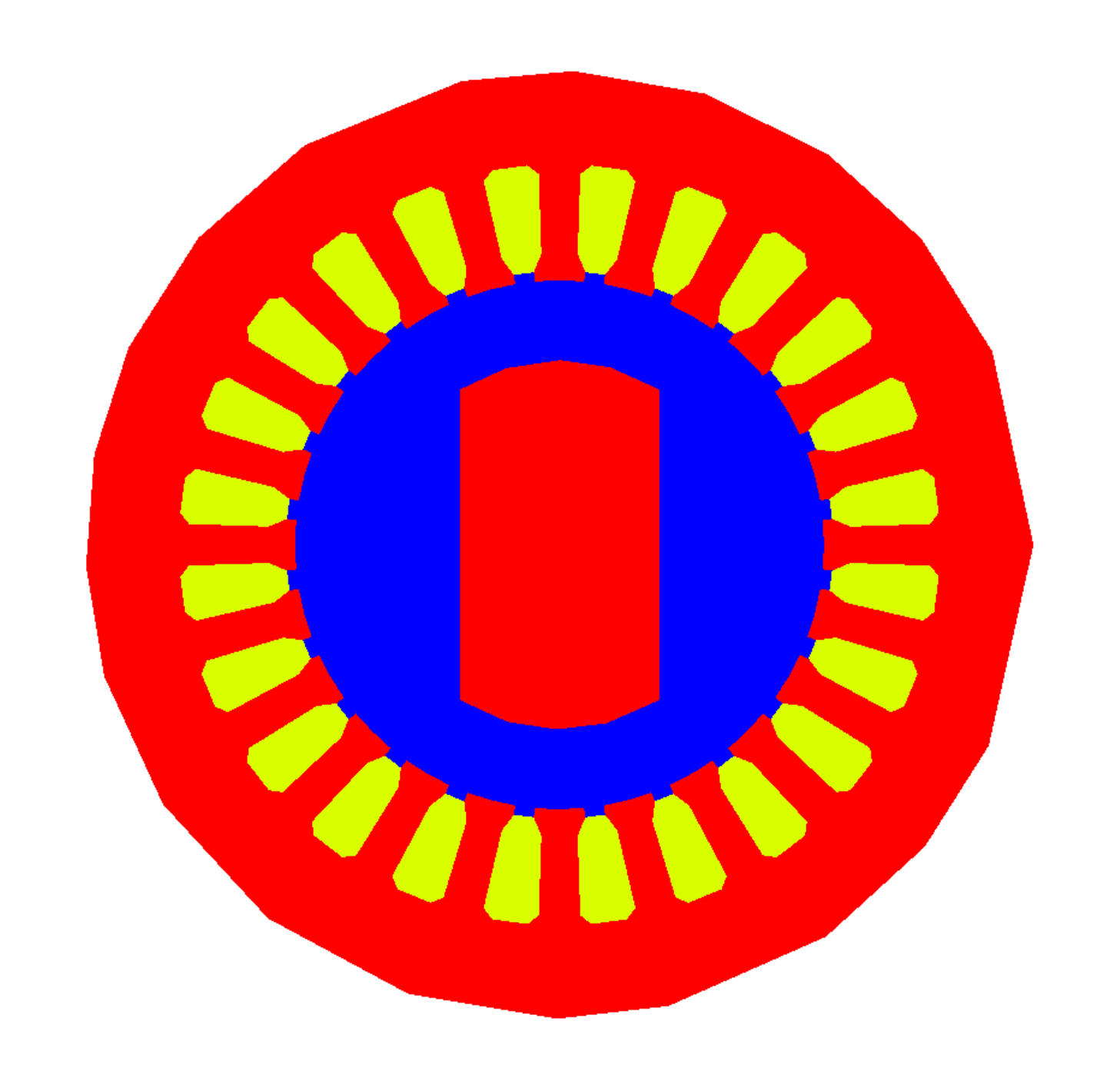}&
    \includegraphics[width=0.55\textwidth]{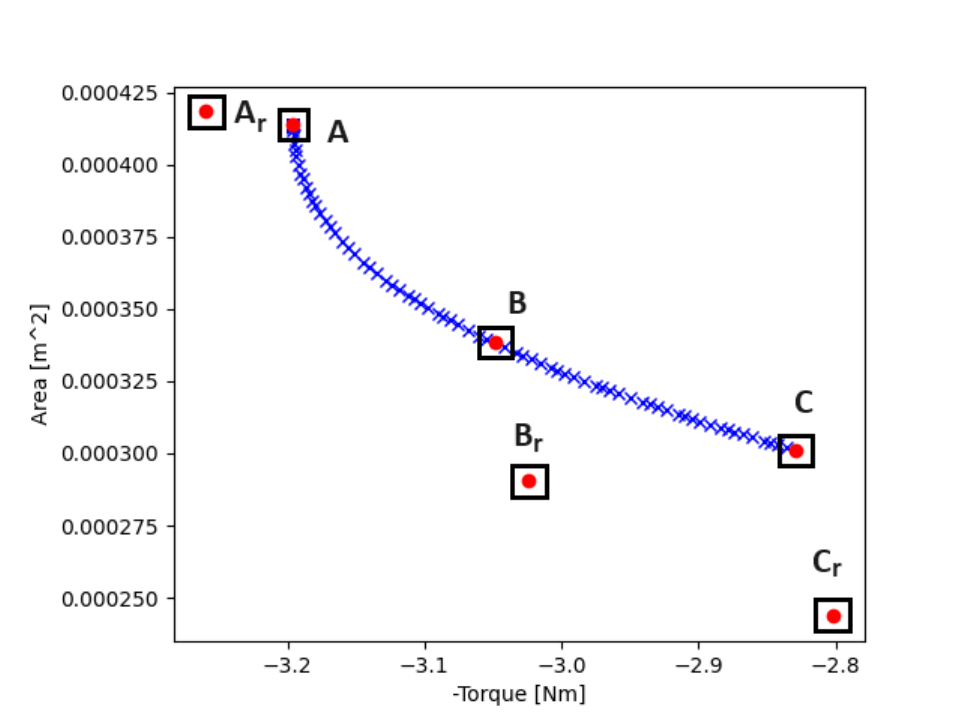}
\end{tabular}
\caption{Left: Initial design. Ferromagnetic material is in red, air is in blue and copper coils are in yellow. Right: Set of Pareto-optimal points traced starting by $A$ in which the negative torque is minimized (homotopy parameter t=0), towards points $B$ (t=0.999779) and $C$ (t=0.999864). Points $A_r$, $B_r$ and $C_r$ are obtained respectively from points $A$, $B$ and $C$ with one uniform refinement of the mesh and following gradient descent optimization.}
\label{cesarano:fig_pareto}
\end{figure}

In this subsection, we present results obtained by applying the procedure detailed in Sect.~\ref{cesarano:sec_paretoshape} to problem \eqref{cesarano:eq_arkkio}--\eqref{cesarano:eq_magnetostatics}.
We start from a simple rotor design constituted by a single piece of ferromagnetic material which is straight on the side and round on the top and the bottom, see Fig.~\ref{cesarano:fig_pareto} (left). We then perform a single objective optimization, maximizing the torque, by a gradient descent algorithm repeatedly using \eqref{cesarano:eq_desc_direction} to obtain the design $\Omega^{(0)}$ solving $\mathcal H(\Omega, 0)(\cdot)=0$. This design optimizes the torque and is \black{referred} to as $A$ in Fig.~\ref{cesarano:fig_pareto} (right). 

\begin{remark}
When using a gradient descent algorithm to obtain the optimized design for the torque, one typically arrives in the proximity of a local minimum but not with very high precision. In our experiments, the residual in the optimality condition was still in the order of $10^{-4}$. Here, one has to hope that this is close enough for the shape Newton to converge and to get back to the homotopy curve. During the homotopy method which is based on Newton updates, the residual is reduced to below $10^{-10}$. 

As an alternative to using a gradient descent method for this initial problem, also here one could apply a homotopy method with a suitably chosen "simple" problem that is solved by the initial design. The choice of such a suitable "simple" problem, however, is in general not clear in the case of shape optimization.
\end{remark}

From here we follow the homotopy defined by \eqref{cesarano:eq_H_motor} by progressively increasing the homotopy parameter $t$ in order to get optimized design which are Pareto-optimal. The increase in the homotopy parameter can be arbitrarily chosen, giving some control on the number of points of the Pareto front and the distance between them, which is desirable from the applications point of view. In general, $t$ can be chosen small enough to get an arbitrary number of points on the Pareto front without large gaps in the objective space. The higher the number of desired points is, the more this approach becomes computationally efficient when compared to multi-start gradient-based methods such as the weighted sum method. The Pareto front points obtained are shown in Fig.~\ref{cesarano:fig_pareto} (right). 

The homotopy method stops at $t = 0.999864$. Due to the different order of magnitude of the 2 considered cost functions, the majority of the Pareto-optimal points are obtained with $t$ very close to 1. Rescaling  the cost functions is possible, if desired. The algorithm stops before actually getting to 1, though. This is expected, because for $t=1$ the problem is to minimize the iron volume, meaning that the rotor should get to only air. This is not possible to obtain without changing the topology. Also, we keep the design fixed at the boundary of the rotor by enforcing deformations there to vanish. The algorithm might go a bit further by allowing the points on the rotor boundary to move tangentially.

In our numerical experiments, we followed the idea of multi-resolution shape optimization \cite{cirak}, meaning to perform first optimization runs on a very coarse mesh. Besides the obvious aspect of efficiency, we also experienced less problems with mesh quality due to large deformations and local minimizers. While the former problem could be remedied by remeshing, the latter aspect was also reported in \cite{cirak}.

\begin{figure}
\includegraphics[width=0.3\textwidth]{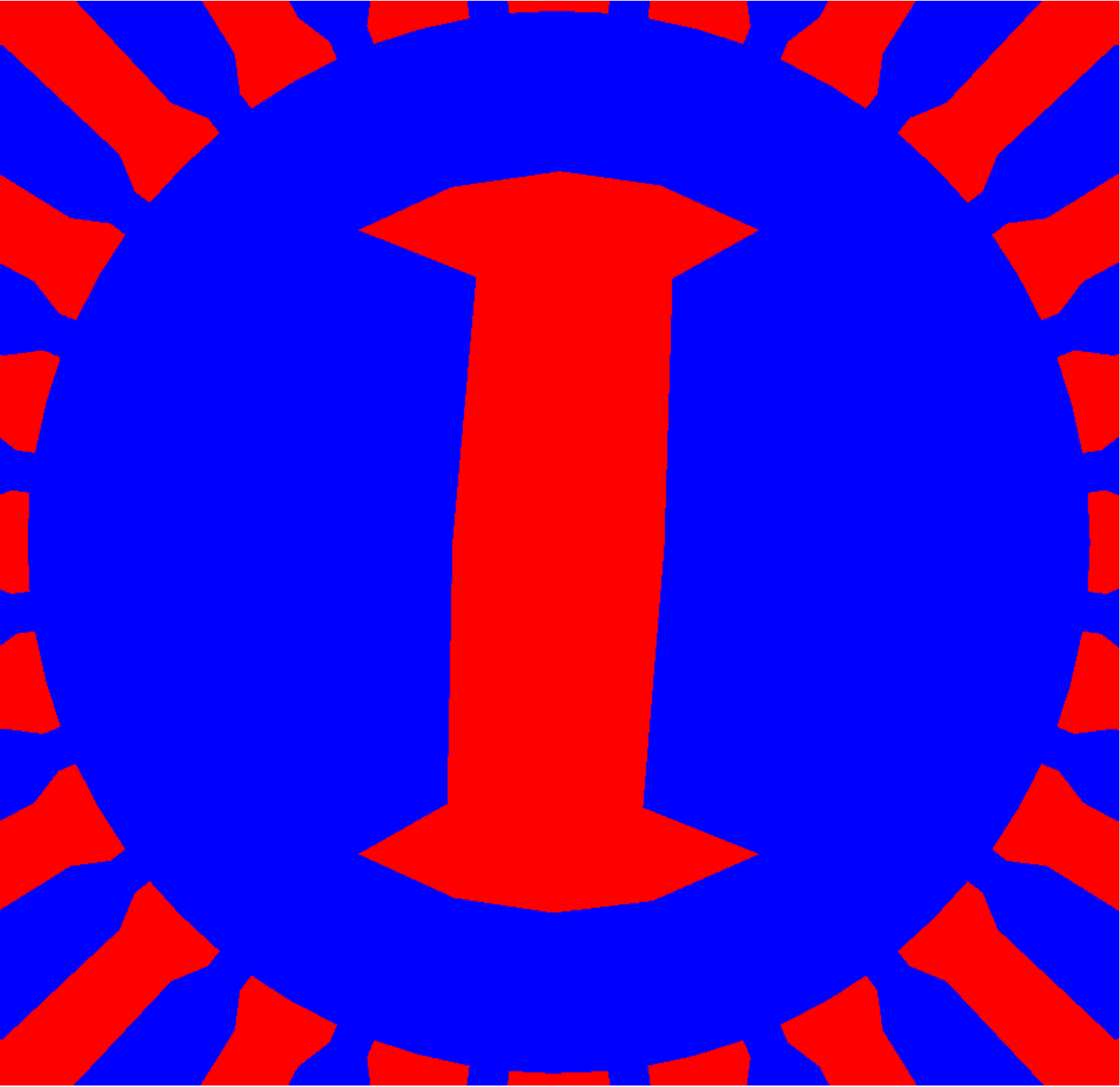}\quad
\includegraphics[width=0.3\textwidth]{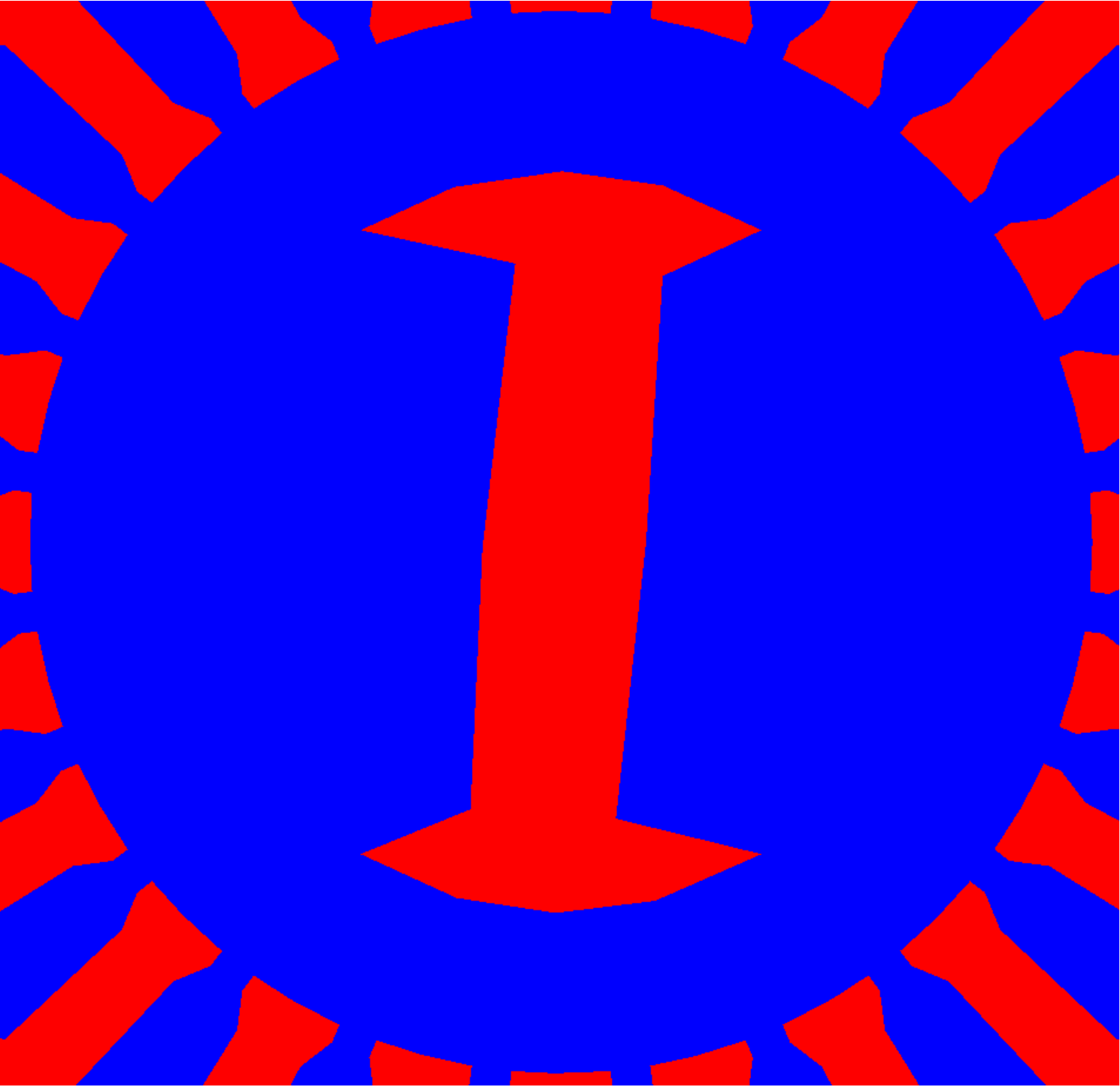}\quad
\includegraphics[width=0.3\textwidth]{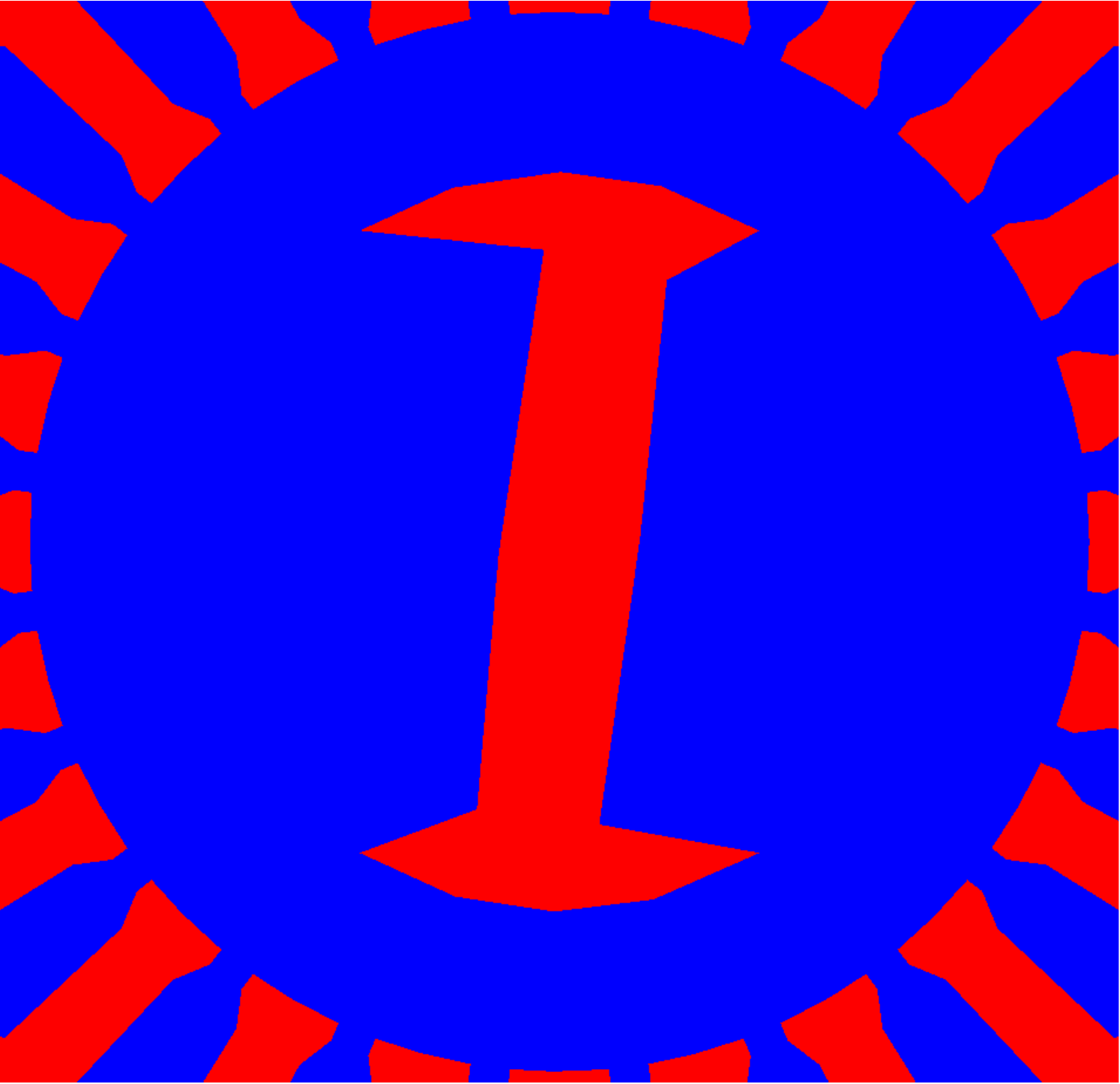}
\caption{Optimized designs corresponding to Pareto-optimal points $A$ (t=0), $B$ (t=0.999779) and $C$ (t=0.999864), from left to right.}
\label{cesarano:fig_optimal}
\end{figure}

\begin{figure}
\includegraphics[width=0.3\textwidth]{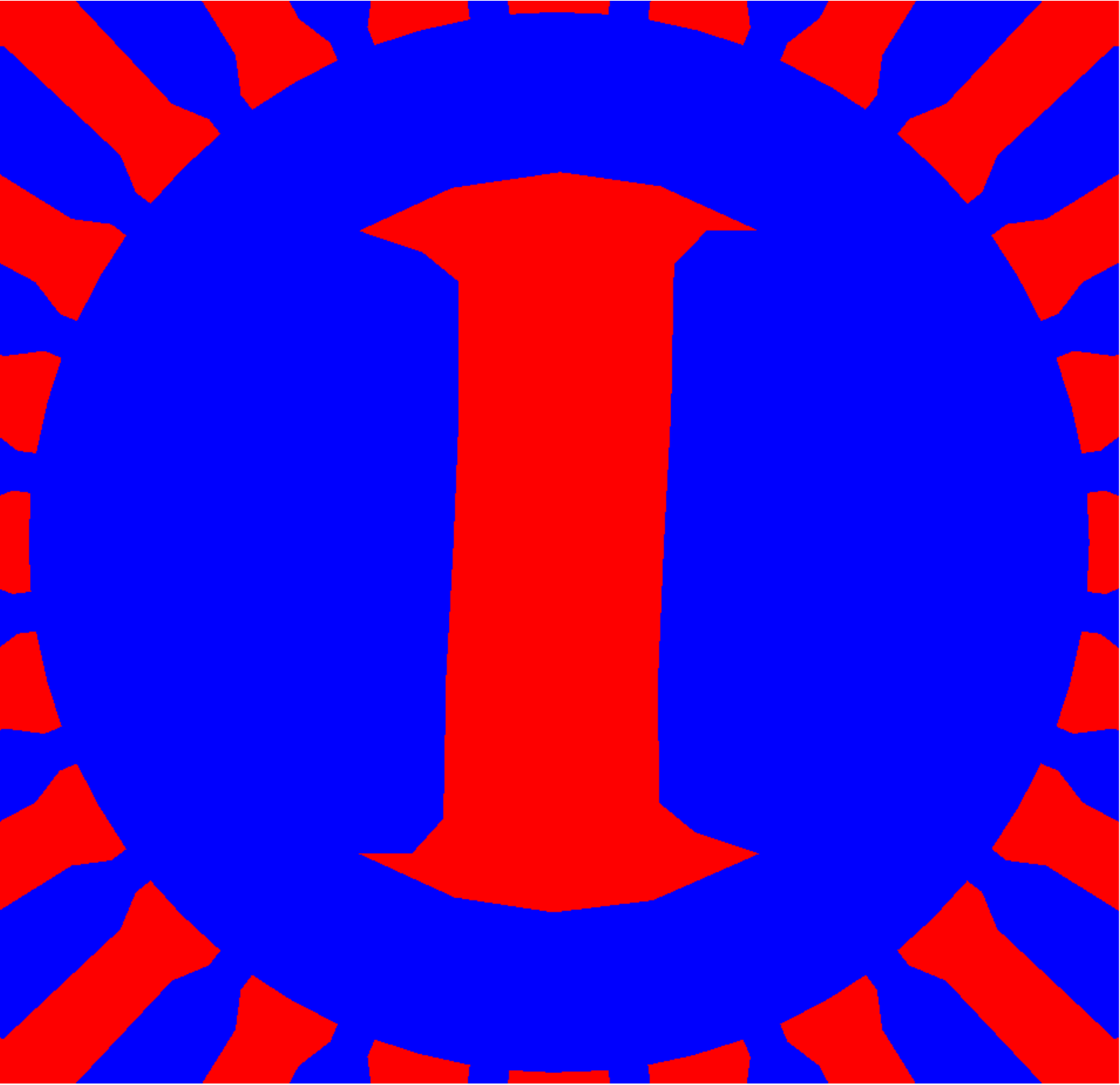}\quad
\includegraphics[width=0.3\textwidth]{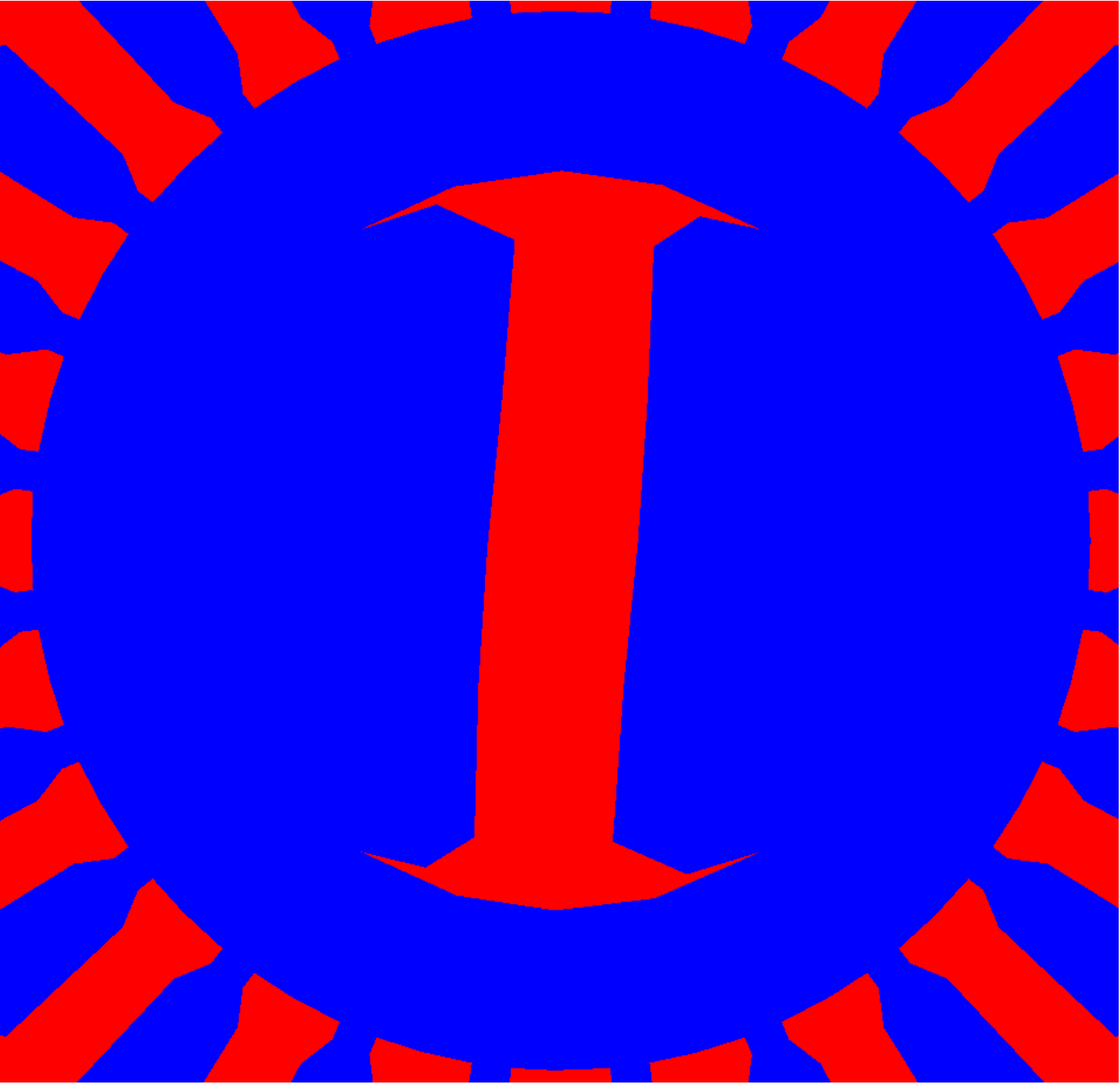}\quad
\includegraphics[width=0.3\textwidth]{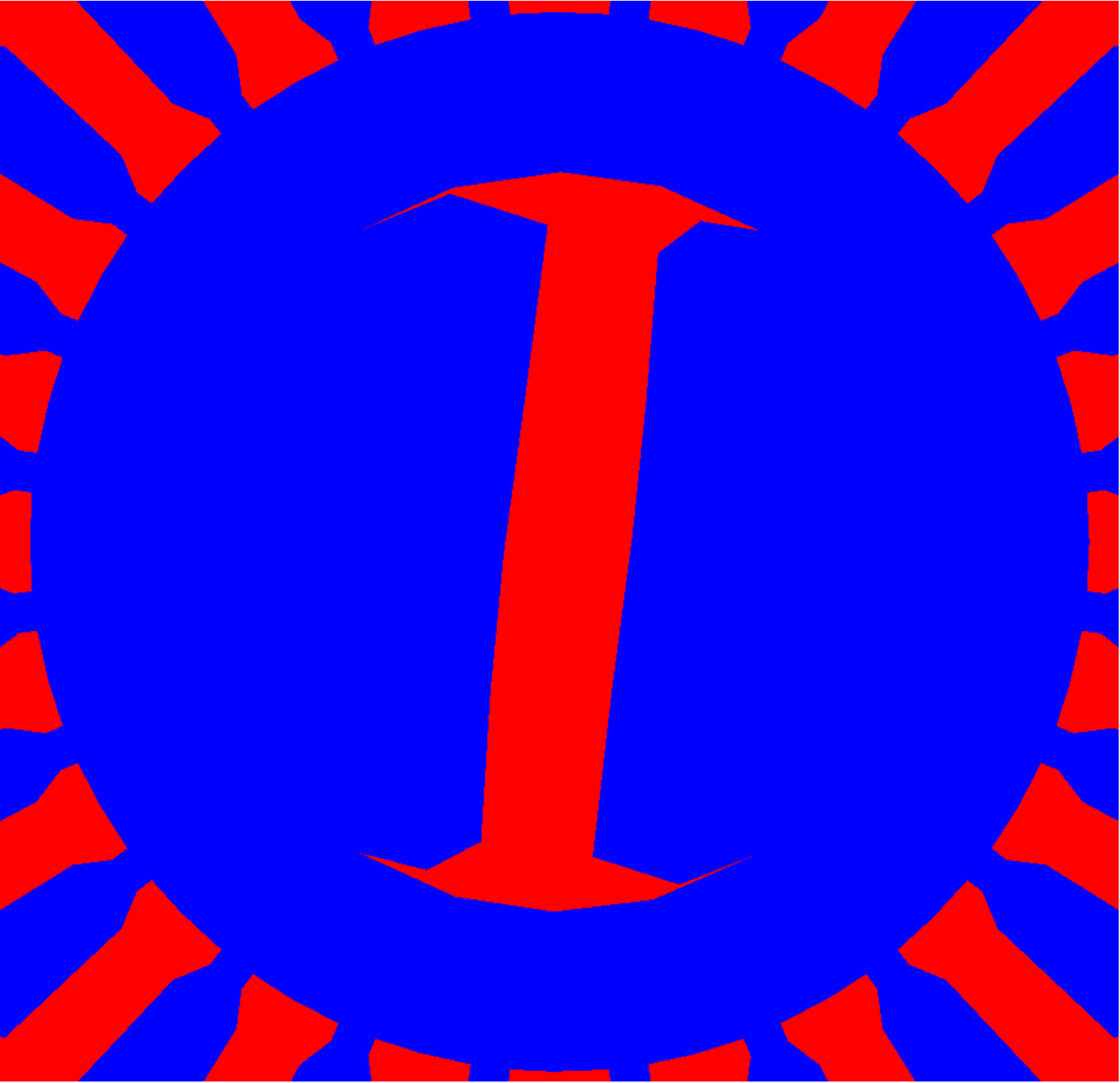}
\caption{Optimized designs corresponding to Pareto-optimal points $A_r$, $B_r$ and $C_r$, from left to right.}
\label{cesarano:fig_optimal_refined}
\end{figure}

In our framework this is translated into finding Pareto-optimal points of a coarse problem, that can be the starting point for a following optimization after a mesh refinement. This is how we obtain the points $A_r$, $B_r$ and $C_r$, whose designs are shown in Fig.~\ref{cesarano:fig_optimal_refined}, starting from the points $A$, $B$ and $C$, shown in Fig.~\ref{cesarano:fig_optimal}. The corresponding objective values are depicted in Fig. \ref{cesarano:fig_pareto} (right).

\section{Conclusions} 

In this work we dealt with multi-objective shape optimization in the context of electric machines and in particular of a synchronous reluctance motor with a single pair of poles. We presented an extension of the Pareto front tracing methods used for gradient descent optimization algorithms to our framework of shape optimization using shape derivatives. In order to do so, we used homotopy and shape Newton methods. We achieved a Pareto front with several Pareto-optimal points, starting from the optimized design, result of one single-objective optimization. We also further optimized 3 selected designs on a finer mesh, in a multi-resolution approach. Further work may include considering other scenarios, with different \black{types} of machines and different cost functions. Also, it is natural to think about the extension of the proposed method to the case of more cost functions, resulting in a Pareto front in more dimensions, with the challenge of the increase in computational cost.

\section*{Acknowledgments}

The work of A. C. and P. G. is supported by the FWF funded project P32911 as well as
the joint DFG/FWF Collaborative Research Centre CREATOR (CRC – TRR361/F90) at TU Darmstadt, TU Graz, RICAM and JKU Linz. 

\bibliographystyle{unsrt}
\bibliography{references}

\end{document}